\documentclass[11pt,reqno,twoside]{article}
\usepackage{amsmath,amssymb,theorem,color,epsfig,epic,subfigure,amsfonts,hhline,graphicx}

\theorembodyfont{\itshape}
\newtheorem{thm}{Theorem}[section]
\newtheorem{lem}[thm]{Lemma}

{\theorembodyfont{\upshape}

\newtheorem{ex}[thm]{Example}
}

\newcommand{\Proof}[1][]{\noindent{\itshape Proof#1. }}
\newcommand{\EndProof}{\hfill$\Box$\bigskip}

\makeatletter
    
    \@addtoreset{equation}{section}
\makeatother

\def\tilde{\widetilde}

    \def\e{\varepsilon}

\begin{document}
\title{Stopping times in the game Rock-Paper-Scissors}
\author{Kyeonghoon Jeong\footnote{Faculty of Liberal Education, Seoul National University, 1 Gwanak-ro, Gwanak-gu, Seoul 08826, Korea. E-mail: zornslemon@snu.ac.kr} and Hyun Jae Yoo\footnote{Department of Applied Mathematics and Institute for Integrated Mathematical Sciences,
Hankyong National University, 327 Jungang-ro, Anseong-si,
Gyeonggi-do 17579, Korea. E-mail: yoohj@hknu.ac.kr}  \footnote{Corresponding author} }
\date{ }
   \maketitle

\begin{abstract}
In this paper we compute the stopping times in the game Rock-Paper-Scissors. By exploiting the recurrence relation we compute the mean values of stopping times. On the other hand, by constructing a transition matrix for a Markov chain associated with the game, we get also the distribution of the stopping times and thereby we compute the mean stopping times again. Then we show that the mean stopping times increase exponentially fast as the number of the participants increases. 
\end{abstract}

\section{Introduction}\label{sec:introduction} 

\begin{figure}[htbp]\label{fig:gawibawibo}
\begin{center}
        \includegraphics[width=0.4\textwidth]{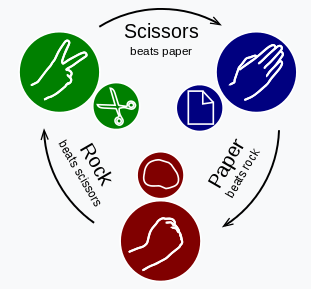}
        \caption{The game Rock-Paper-Scissors ({\small picture taken from Wikipedia}) }
  \end{center}
\end{figure}

The game Rock-Paper-Scissors is perhaps the most famous game known world widely for choosing a winner among participants. The rule is very simple. Each player shows by hand one of Rock, Paper, or Scissors at the same time. The Rock beats the Scissors and loses to the Paper. The Scissors  beat the Paper and lose to Rock, and the Paper beats the Rock and loses to the Scissors. See the picture in Figure \ref{fig:gawibawibo}. It is natural to ask the ending time of the game when there are a certain  number of participants in the game. 

The purpose of this paper is to answer this question. When a number of participants is fixed, we will discuss the following questions: (i) the distribution of ending time, (ii) the mean ending time. We will also investigate the asymptotic behavior of the mean ending time as the number of participants increases. To summarize the results, we have obtained a concrete formula for the distribution and the mean value of the ending time and showed that the mean ending time increases exponentially fast as the number of participants increases. 

This paper is organized as follows. In Section 2, we exploit the recurrence relation for the mean stopping times of the game. Using exponential generating function, we compute the mean stopping times. 
In Section 3 we introduce a Markov chain for the game. The transition matrix of this Markov chain will have crucial roles in the computations of our interests. Using the transition matrix, we represent the mass functions of the stopping times and compute the mean stopping times. In Section 4, we discuss the asymptotic behavior of the mean stopping  times. In Section 5, we give the proofs for the main results.     

\section{Recurrence relation}\label{sec:recurrence_relation}

In this section we compute the mean ending time of the game started with a certain number of participants. Here we will elucidate the recurrence relation for the stopping times. In the next section we will introduce  a different method.

For $n\ge 2$, let $\tau_n$ be the ending time of the game with $n$ participants, i.e., the time of the final winner is determined. Let $E_n:=\mathbb E(\tau_n)$ be the expectation of $\tau_n$. For simplicity we let $\tau_1=0$, i.e., $E_1=0$. For $n\ge 2$ and $j=1,\cdots,n$, let $p(n,j)$ be the probability that $j$ persons survive after one round of the game with $n$ participants. We then have the relation: for $n\ge 2$ and $k\ge 1$
\begin{equation}\label{eq:prob_stopping}
\mathbb P(\tau_n=k)=\sum_{j=1}^np(n,j)\mathbb P(\tau_j=k-1).
\end{equation}
We thus get the recurrence relation for the expectation values:
\begin{equation}\label{eq:recurrence_relation}
  E_n =\displaystyle \sum_{j=1}^{n} p(n,j) (E_j + 1), \quad n\ge 2.
\end{equation}
From the rule of the game we easily see that 
\begin{equation}\label{eq:reducing_probability}
p(n,j)=\begin{cases}\frac1{3^{n-1}}{n \choose j},&j<n\\ 1-2\left(\frac23\right)^{n-1}\left(1-\frac1{2^{n-1}}\right),&j=n
\end{cases}.
\end{equation}
So we obtain
\begin{align}\label{recurrence}
(2^n - 1) E_n = 3^{n-1} + \displaystyle \sum_{j=1}^n \binom n j   E_j,\quad n\ge 2. 
\end{align} 
From this recurrence relation with the initial condition $  E_1 = 0$,
we can easily calculate
\[  E_2 = \frac 3 2, \   E_3 = \frac 9 4 , \   E_4 = \frac {45}{14}, \   E_5 = \frac {157}{35} , \]
etc. 

Let us now compute the general formula for $E_n$ solving the equation \eqref{recurrence}. 
Put $E_0 = 0$ for convenience. Define the exponential generating function of $\{E_n\}$ by
\begin{equation}\label{eq:exp_gen_function}
E(x) := \sum_{n=0}^\infty \frac {  E_n}{n!} x^n=\sum_{n=2}^\infty \frac {  E_n}{n!} x^n.
\end{equation}  
By \eqref{recurrence} it satisfies the following functional equation
\[ E(2x) - E(x) = \frac 1 3 (e^{3x} - 1 - 3x) + e^x E(x) , \]
that is,
\[ E(2x) = (e^x + 1) E(x) + \frac 1 3 (e^{3x} -1 - 3x) . \]
Now we have
\[ \frac {E(2x)}{e^{2x} - 1} = \frac {E(x)}{e^{x} - 1} + \frac 1 3 \frac {e^{3x} -1-3x}{e^{2x} - 1}. \]
Let $h(x) = \dfrac 1 3 \dfrac {e^{3x} - 1 - 3x}{e^{2x} - 1} = \displaystyle \sum_{n=1}^\infty \frac {h_n}{n!} x^n$. By calculation, we have $h_1 = \frac 3 4 , h_2 = 0, h_3 = \frac 3 8 , h_4 = \frac 3 5, h_5 = \frac 1 4, h_6 = - \frac 3 7$, etc. (And $h_n$ can be expressed using Bernoulli numbers.)

If we put $F(x) = \dfrac {E(x)}{e^{x} - 1}$, then $F(2x) = F(x) + h(x)$ implies
\[ 2^n F^{(n)} (0) = F^{(n)} (0) + h_n .\]
Therefore we have
\begin{equation}\label{eq:exp_gen_function_sol}
 \frac {E(x)}{e^{x} - 1} = \sum_{n=1}^\infty \frac {h_n}{n! (2^n - 1)} x^n . \end{equation}
Or equivalently $F(x) = \displaystyle \sum_{k=1}^\infty h \left( \frac x {2^k} \right)$. Finally, expanding \eqref{eq:exp_gen_function_sol}, we arrive at 
\begin{thm}\label{thm:mean_stopping_time_via_recurrence} In the game Rock-Paper-Scissors with $n$ participants, $n\ge 2$, the mean stopping time is given by 
\begin{align}\label{eq:mean_stopping_time}
  E_n = \sum_{k=1}^{n-1} \binom n k \frac {h_k}{2^k - 1}.
\end{align}
\end{thm}
For example, we have $E_4 = 4 \frac {h_1} 1 + 6 \frac {h_2} 3 + 4 \frac {h_3} 7 = 3 + \frac 3 {14}$ once again.


\section{Markov chain of the game Rock-Paper-Scissors and stopping times}

In this section we investigate the stopping times $\tau_n$ with a different method. It turns out to be convenient to introduce a 
  stochastic matrix $P$ in the state space $\mathbb N$, the set of natural numbers, as follows.
\[
P:=[p(i,j)]_{i,j\in \mathbb N},
\]
where, as before, 
$$
p(i,j):=\mathbb{P}(j \text{ persons survive after one round of the game of } i \text{ persons}),
$$ 
and it was computed in \eqref{eq:reducing_probability}.
Here we take $p(i,j)=0$ for $j>i$ and we put $p(1,1)\equiv 1$.  
Some of the leading terms of $P$ are as follows.
\[
P=\left[\begin{matrix}1&0&0&0&0&\cdots\\
\frac23&\frac13&0&0&0&\cdots\\
\frac13&\frac13&\frac13&0&0&\cdots\\
\frac{4}{27}&\frac{2}{9}&\frac{4}{27}&\frac{13}{27}&0&\cdots\\
\vdots&\vdots&\vdots&\vdots&\vdots&\ddots\end{matrix}\right].
\]
It is obvious that for the Markov chain with the transition matrix $P$, the state $1$ is recurrent (actually absorbing) and all the other states are transient \cite{N}. Thus in the Markov chain there is only one trivial invariant measure $\delta_1$, the Dirac measure at the point $1$. We are interested in the ending time of the game. This is the hitting time of the state $1$.  

Let $\{{\bf e}_i\}_{i\in \mathbb N}$ be the canonical basis of the Hilbert space $(\mathcal H,\langle \cdot,\cdot\rangle):=l^2(\mathbb N)$. For $n\ge 2$, let $S_{[2,n]}$ be the projection onto the subspace  $\mathrm{span}\{{\bf e}_i:2\le i\le n\}$ and we define  
\[
P_{n}:=S_{[2,n]}PS_{[2,n]}.
\]
Let $p(\cdot,1):=(p(i,1))_{i\in \mathbb N}\in \mathcal H$ and denote $\psi_n:=S_{[2,n]}p(\cdot,1)$. Let us define the mass function of the distribution of $\tau_n$:
\[
p_n(k):=\mathbb{P}(\tau_n=k), \quad k\in \mathbb{N}.
\] 
We have the following results.
\begin{thm}\label{thm:probabilities}
In the game Rock-Paper-Scissors, let $\tau_n$ be the stopping time of ending of the game. Then, 
\begin{enumerate}
\item[(i)] $\tau_n$ has the mass function
\[
p_n(k)=\langle {\bf e}_n,P_{n}^{k-1}\psi_n\rangle=P_{n}^{k-1}\psi_n(n).
\]\\[-3ex]
\item[(ii)] The moment generating function of $\tau_n$ is
\[
M_n(t):=\mathbb{E}(e^{t\tau_n})=\langle {\bf e}_n,\frac{e^t}{I-e^tP_n}\psi_n\rangle=\frac{e^t}{I-e^tP_n}\psi_n(n),
\]
and hence the expectation and variance of $\tau_n$ are respectively given by 
\begin{eqnarray*}
E_n=\mathbb{E}(\tau_n)&=&\langle {\bf e}_n,(I-P_{n})^{-2}\psi_n\rangle=\frac1{(I-P_n)^{2}}\psi_n(n)\\
\mathrm{Var}(\tau_n)&=& \langle {\bf e}_n,\frac{I+P_n}{(I-P_n)^3}\psi_n\rangle-\langle {\bf e}_n,(I-P_{n})^{-2}\psi_n\rangle^2\\
&=&\frac{I+P_n}{(I-P_n)^3}\psi_n(n)-\left(\frac1{(I-P_n)^{2}}\psi_n(n)\right)^2.
\end{eqnarray*}
\end{enumerate}
\end{thm}
\begin{ex} We compute the mean ending times for some values of $n$. For $n=2,3,4$, 
\begin{eqnarray*}
P_{2}&=&\left[\begin{matrix}\frac13\end{matrix}\right], \quad \psi_2=\left[\begin{matrix}\frac23\end{matrix}\right], \\
P_{3}&=&\left[\begin{matrix}\frac13&0\\\frac13&\frac13\end{matrix}\right],\quad \psi_3=\left[\begin{matrix}\frac23\\\frac13\end{matrix}\right]\\
P_4&=&\left[\begin{matrix}\frac13&0&0\\\frac13&\frac13&0\\\frac{2}{9}&\frac{4}{27}&\frac{13}{27}\end{matrix}\right], \quad \psi_4=\left[\begin{matrix}\frac23\\\frac13\\\frac{4}{27}\end{matrix}\right].
\end{eqnarray*}
Therefore,
\[
E_2=\frac1{\left(1-\frac13\right)^2}\left(\frac23\right)=\frac32,
\]
\[
E_3=\frac1{(I-P_3)^{2}}\psi_3(3)=\frac94\left[\begin{matrix}1&0\\1&1\end{matrix}\right]
\left[\begin{matrix}\frac23\\\frac13\end{matrix}\right](3)=\frac94,
\]
\[
E_4=\frac1{(I-P_4)^{2}}\psi_4(4)=\left[\begin{matrix}\frac94&0&0\\\frac94&\frac94&0\\
\frac{639}{196}&\frac{72}{49}&\frac{729}{196}\end{matrix}\right]\left[\begin{matrix}\frac23\\\frac13\\\frac{4}{27}\end{matrix}\right](4)=\frac{45}{14}.
\]
These are the same as the corresponding values computed in Section 2. 
\end{ex}
 
\section{Asymptotic behavior}

In this section we discuss the asymptotic behavior of the mean value of ending time. The following is the main result for the asymptotics.
\begin{thm}\label{thm:asymptotics}
In the game Rock-Paper-Scissors, the mean ending time increases exponentially fast as the number of participants increases. More precisely, we have
\[
\frac13\left(\frac32\right)^n\le\mathbb{E}(\tau_n) \le \frac13n^3\left(\frac32\right)^n.
\]
\end{thm}
The lower and upper bounds we will get by considering the first exit times from the present number of participants in the game. For the exponential growth, however, we will introduce also a different method  by which we can learn a little bit how the dynamics of the game proceeds.\footnote{After the paper has been accepted for publication, by a discussion with Professor Tomoyuki Shirai, it was known that the asymptotic behavior could be substantially improved. See the Note Added in Proof.}

\subsection{Lower bound}\label{subsec:lower_bound}

\subsubsection{Trajectory of the dynamics}

First we notice that since the mean ending time is represented by the inverses of triangular matrices, it is very helpful to have knowledge of them.\footnote{We have taken Lemma \ref{lem:inverse_triangular_matrix} and its proof from a note posted by Robert Lewis in google search.}
\begin{lem}\label{lem:inverse_triangular_matrix}
Let $A$ be an $n\times n$ triangular matrix with non-zero diagonal elements. If we put $A=D+R$, where $D$ is the diagonal part of $A$, then we have
\[
A^{-1}=\sum_{k=0}^{n-1}(-D^{-1}R)^kD^{-1}.
\]
\end{lem}
\Proof
We can write $A=D(I+D^{-1}R)$. We notice that $D^{-1}R$ is strictly triangular. Furthermore, it is nilpotent, i.e., $(D^{-1}R)^n=0$. Now we have $A^{-1}=(I+D^{-1}R)^{-1}D^{-1}$. By using the identity
\[
(1+x)\sum_{k=0}^m(-x)^k=1-(-x)^{m+1}
\]
for $x=D^{-1}R$ and $m=n-1$, and the nilpotency of $D^{-1}R$, we have
\[
(I+D^{-1}R)^{-1}=\sum_{k=0}^{n-1}(-D^{-1}R)^k.
\]
This completes the proof.
\EndProof\\

From now, let us show the exponential growth of the mean ending time. 
By Theorem \ref{thm:probabilities} we have 
\begin{eqnarray}\label{eq:expectation_with_matrix_components}
\mathbb{E}(\tau_n)&=&\langle {\bf e}_n,(I-P_{n})^{-2}\psi_n\rangle\nonumber\\
&=&\sum_{j=2}^n(I-P_{n})^{-2}(n,j)\psi_n(j).
\end{eqnarray}
Notice that the matrix $I-P_n$ is a lower triangular matrix of size $(n-1)\times (n-1)$. We can write
\[
I-P_n=D_n-L_n,
\]
where $D_n$ is the diagonal part of $I-P_n$. We notice here that $D_n$ has strictly positive components and $L_n$ is strictly lower triangular with nonnegative components. By Lemma \ref{lem:inverse_triangular_matrix} we have 
\begin{eqnarray*}
(I-P_n)^{-1}&=&\sum_{k=0}^{n-2}(D_n^{-1}L_n)^kD_n^{-1}\\
&=&D_n^{-1/2}\left[\sum_{k=0}^{n-2}\left(D_n^{-1/2}L_nD_n^{-1/2}
\right)^k\right]D_n^{-1/2}.
\end{eqnarray*}
We thus get
\begin{equation}\label{eq:expansion}
(I-P_n)^{-2}=D_n^{-1/2}\left[\sum_{k=0}^{n-2}\left(D_n^{-1/2}L_nD_n^{-1/2}
\right)^k\right]   D_n^{-1}\left[\sum_{l=0}^{n-2}\left(D_n^{-1/2}L_nD_n^{-1/2}
\right)^l\right]D_n^{-1/2}.
\end{equation}
Since all the matrix components of $D_n$ and $L_n$ are nonnegative and also $\psi_n$ has positive components, we see from equations \eqref{eq:expectation_with_matrix_components} and \eqref{eq:expansion} that (taking just a single term with $j=[n/2]$ in \eqref{eq:expectation_with_matrix_components}, and $k=0$ and $l=1$ in \eqref{eq:expansion})
\begin{equation}\label{eq:lower_bound}
\mathbb{E}(\tau_n) \ge \left(D_n^{-2}L_nD_n^{-1}\right)(n,[n/2])\psi_n([n/2]),
\end{equation}
where $[r]$ means the integer part of a real number $r$. For $2\le j\le n$, the $j$th component of $D_n$, denoted by $D_n(j)$, is given by 
\[
D_n(j)=1-p(j,j)=2\left(\frac23\right)^{j-1}\left(1-\frac1{2^{j-1}}\right),
\]
and 
\[
\psi_n(j)=p(j,1)=\frac{j}{3^{j-1}}.
\]
For $2\le i\le n$ and $2\le j<i$, we have 
\[
L_n(i,j)=P_n(i,j)=\frac1{3^{i-1}}{i\choose j}.
\]
Therefore, plugging into \eqref{eq:lower_bound} we get
\[
\mathbb{E}(\tau_n) \ge 2^{-2}\left(\frac32\right)^{2(n-1)}\frac1{3^{n-1}}{n\choose {[n/2]}}\,2^{-1}\left(\frac32\right)^{[n/2]-1}
\frac{[n/2]}{3^{[n/2]-1}}.
\]
By Stirling's formula we have
\[
{n\choose {[n/2]}}\sim \frac{2^{n+1}}{\sqrt{2\pi n}}.
\]
Thus finally we get 
\[
\mathbb{E}(\tau_n) \ge C\sqrt{n}\left(\frac3{2^{3/2}}\right)^n,
\]
where $C$ is a constant. This shows that the mean ending time grows exponentially as the number of participants increases.

\subsubsection{Exit times}

For each $n\ge 2$, let $T_{\mathrm{ex}}^{(n)}$ be the first exit time from the initial state in the game starting with $n$ participants, namely,  
\begin{equation}\label{eq:exit_time}
T_{\mathrm{ex}}^{(n)}:=\inf\{k\ge 1: \#\text{of participants}<n\text{ after }k\text{th round of the game}\}.
\end{equation}
We have
\begin{equation}\label{eq:exit_probability}
\mathbb{P}(T_{\mathrm{ex}}^{(n)}=k)=\sum_{j=1}^{n-1}p(n,n)^{k-1}p(n,j)=p(n,n)^{k-1}-p(n,n)^k.
\end{equation}
Obviously we have
\begin{eqnarray*}
\mathbb{E}(\tau_n)&\ge &\mathbb{E}(T_{\mathrm{ex}}^{(n)})\\
&=&\sum_{k=1}^\infty k\mathbb{P}(T_{\mathrm{ex}}^{(n)}=k)\\
&=&\sum_{k=1}^\infty k\left(p(n,n)^{k-1}-p(n,n)^k\right)\\
&=&\frac1{1-p(n,n)}\ge \frac13\left(\frac32\right)^n.
\end{eqnarray*}
This proves the lower bound.

\subsection{Upper bound}\label{subsec:upper_bound}

In this subsection we compute the upper bound for the asymptotic mean stopping times of the game. The basic idea is to look at carefully the trajectory of decreasing numbers of participants as the game goes on. 

Suppose that the game starts with $n$ participants, call it an $n$-block game. As we have seen in the former subsection, we have to wait a certain time, say $n_0\ge 0$, until the game firstly exits the $n$-block, then it goes into a small size, say $j_1$-block game. There we wait another exit time, say $n_1$, and then the game goes into further smaller block. The game continues this way and at a certain time it at last jumps into 1-block, the end point. The number of jumps into smaller sized blocks runs between 1 and $n-1$. Therefore, we can compute the mean stopping time of the game as follows: (below we denote  $\mathbb N_0:=\{0\}\cup\mathbb N)$
\begin{eqnarray} \label{eq:expectation_stopping_time}
&&\mathbb E(\tau_n)\\
&=&\sum_{k=1}^{n-1}\sum_{1<j_{k-1}<\cdots<j_1<n}\sum_{(n_0,\cdots,n_{k-1})\in \mathbb N_0^k}(k+n_0+\cdots+n_{k-1})p(n,n)^{n_0}\cdots p(j_{k-1},j_{k-1})^{n_{k-1}}\nonumber \\ &&\hskip 3true cm \times p(n,j_1)p(j_1,j_2)\cdots p(j_{k-1},1).\nonumber
\end{eqnarray}
Taking a change of variables $n_j+1\to n_j$ for $j=0,\cdots, k-1$, and using the formula $p(i,j)$ in \eqref{eq:reducing_probability}, particularly $p(i,i)=1-2\left(\frac23\right)^{i-1}\left(1-\frac1{2^{i-1}}\right)\le 1-\left(\frac23\right)^{i-1}$, we get
\begin{eqnarray*}  
&&\mathbb E(\tau_n)\\
&\le&\sum_{k=1}^{n-1}\sum_{1<j_{k-1}<\cdots<j_1<n}\sum_{(n_0,\cdots,n_{k-1})\in \mathbb N^k}(n_0+\cdots+n_{k-1})\left(1-\left(\frac23\right)^{n-1}\right)^{n_0-1}\cdots\\
&&\hskip 1 true cm \times\left(1-\left(\frac23\right)^{j_{k-1}-1}\right)^{n_{k-1}-1}
 \frac1{3^{n-1}}\frac1{3^{j_1-1}}\cdots\frac1{3^{j_{k-1}-1}}{n\choose j_1}{j_1\choose j_2}\cdots{j_{k-1}\choose 1}\\
 &=&\sum_{k=1}^{n-1}\sum_{1<j_{k-1}<\cdots<j_1<n}\left[\left(\frac32\right)^{n-1}
 +\cdots+\left(\frac32\right)^{j_{k-1}-1}\right]\frac1{2^{n-1}}\frac1{2^{j_1-1}}
 \cdots\frac1{2^{j_{k-1}-1}}\\
&&\hskip 1 true cm \times {n\choose j_1}{j_1\choose j_2}\cdots{j_{k-1}\choose 1}\\
&\le& \left(\frac32\right)^{n-1} \sum_{k=1}^{n-1}k\sum_{1<j_{k-1}<\cdots<j_1<n}\frac1{2^{n-1}}\frac1{2^{j_1-1}}
 \cdots\frac1{2^{j_{k-1}-1}} {n\choose j_1}{j_1\choose j_2}\cdots{j_{k-1}\choose 1}.
\end{eqnarray*}
We see that 
\[
{n\choose j_1}{j_1\choose j_2}\cdots{j_{k-1}\choose 1}=\frac{n!}{(n-j_1)!(j_1-j_2)!\cdots(j_{k-1}-1)!}.
\]
Using this and taking a change of variables $j_i-1\to j_i$ for $i=1,\cdots,k-1$, and applying multinomial expansion, we get
\begin{eqnarray*}  
&&\mathbb E(\tau_n)\\
&\le&n\left(\frac32\right)^{n-1} \sum_{k=1}^{n-1}k\sum_{0<j_{k-1}<\cdots<j_1<n-1}\frac{(n-1)!}{(n-1-j_1)!(j_1-j_2)!\cdots(j_{k-1})!}\\
&&\hskip 1 true cm \times \left(\frac12\right)^{n-1-j_1}\left(\frac1{2^2}\right)^{j_1-j_2}\cdots\left(\frac1{2^k}\right)^{j_{k-1}}\\
&\le& n\left(\frac32\right)^{n-1} \sum_{k=1}^{n-1}k\left(\frac12+\frac1{2^2}+\cdots
+\frac1{2^k}\right)^{n-1}\\
&\le&\frac13n^3\left(\frac32\right)^n.
\end{eqnarray*}
We now have shown Theorem \ref{thm:asymptotics}.
\section{Proof of Theorem \ref{thm:probabilities}}

In this section, we provide with the proofs of the main results.\\
\Proof[ of Theorem \ref{thm:probabilities}] (i) For $n\ge 2$ and $k\ge 1$, it holds that
\begin{eqnarray*}
p_n(1)&=&p(n,1)=\frac{3n}{3^n}\\
p_n(k)&=&\sum_{j_1=2}^n p(n,j_1)p_{j_1}(k-1)\\
&=&\sum_{j_1=2}^n\sum_{j_2=2}^{j_1}p(n,j_1)p(j_1,j_2)p_{j_2}(k-2)=\sum_{j_1=2}^n\sum_{j_2=2}^{ n}p(n,j_1)p(j_1,j_2)p_{j_2}(k-2)\\
&=&P_n^2p_{\cdot}(k-2).
\end{eqnarray*}
In the third equation we have used the fact that $p(j_1,j_2)=0$ for $j_2>j_1$. Repeating the argument we get the result.\\
(ii) We use the above result to compute the moment generating function of $\tau_n$, $M_n(t):=\mathbb{E}(e^{t\tau_n})$. 
\begin{eqnarray*}
M_n(t)&=&\sum_{k=1}^\infty e^{tk}p_n(k)\\
&=&\sum_{k=1}^\infty e^{tk}\langle {\bf e}_n,P_n^{k-1}\psi_n\rangle\\
&=& e^t\sum_{k=1}^\infty\langle {\bf e}_n,\left(e^tP_n\right)^{k-1}\psi_n\rangle\\
&=&\langle {\bf e}_n,\frac{e^t}{I-e^tP_n}\psi_n\rangle.
\end{eqnarray*}
Notice that since the eigenvalues of $P_n$ lie in the open interval $(0,1)$, $M_n(t)$ is well defined in the neighborhood of $t=0$. The mean value and the variance of $\tau_n$ can be computed by differentiating the function $M_n(t)$. The proof is completed.
\EndProof \\
{\bf Note Added in Proof}\\
Here we give an improved result for the asymptotic behavior. We are grateful to  Professor Tomoyuki Shirai for giving us the comments and idea for the improvement.  

For two sequences $f(n)$ and $g(n)$ we write $f(n)\sim g(n)$, as usual, meaning that $\lim_{n\to \infty}f(n)/g(n)=1$.
\begin{thm}\label{thm:asymptotics_improved}
Let $E_n:=\mathbb{E}(\tau_n)$ be the mean ending time of the game Rock-Paper-Scissors started with $n$ participants. Then we have 
\[
E_n=\frac13\left(\frac32\right)^n+r_n,
\]
where 
\[
r_n=\frac13\frac1{2^n-1}\left(\left(\frac32\right)^n+\sum_{s=2}^n
{n\choose s}\left(\frac32\right)^s\sum_{l=1}^\infty l^{-s}2^{\delta(l)n}\right),
\]
with $\delta(l)$ the fractional part of $\log_2l$: $0\le \delta(l)=\log_2l-[ \log_2]<1$. The remainder $r_n$ satisfies $ r_n=\mathrm{o}\left((3/2)^n\right)$, and hence particularly  $\mathbb{E}(\tau_n)\sim\frac13\left(\frac32\right)^n$, with the latter the lower bound.
\end{thm}
\Proof
We recall some functions and their properties discussed in Section \ref{sec:recurrence_relation}. 
\begin{equation}\label{eq:function_h}
h(x) = \dfrac 1 3 \dfrac {e^{3x} - 1 - 3x}{e^{2x} - 1} = \displaystyle \sum_{n=1}^\infty \frac {h_n}{n!} x^n,
\end{equation}
\begin{equation}\label{eq:exp_gen_function_sol_1}
F(x):=\frac {E(x)}{e^{x} - 1} = \sum_{n=1}^\infty \frac {h_n}{n! (2^n - 1)} x^n.
\end{equation}
We have shown the relation 
\begin{equation}\label{eq:expansion_of_Fr}
F(x) = \displaystyle \sum_{k=1}^\infty h \left( \frac x {2^k} \right).
\end{equation} 
Let us introduce the function $\tilde E(x)$ by the formula:
\begin{equation}\label{eq:tilde_E}
\frac{\tilde E(x)}{e^x-1}=\sum_{n=1}^\infty\frac{h_n}{n!2^n}x^n=\frac 1 3 \dfrac {e^{\frac32x} - 1 - \frac32x}{e^{x} - 1}.
\end{equation}
Therefore, we get
\begin{equation}\label{eq:function_tilde_E_n}
\tilde{E}(x)= (e^x-1)h\left(\frac{x}2\right)=\frac 1 3 \left({e^{\frac32x} - 1 - \frac32x}\right) =\frac13\sum_{n=2}^\infty\frac1{n!}\left(\frac32\right)^nx^n.
\end{equation}
Thus, $\tilde E(x)$ is the exponential generating function of the sequence $\tilde E_n:=\frac13\left(\frac32\right)^n$.
On the other hand, from \eqref{eq:exp_gen_function_sol_1} and \eqref{eq:tilde_E} we have 
\[
\frac{E(x)-\tilde{E}(x)}{e^x-1}= \sum_{n=1}^\infty \frac {h_n}{n! 2^n(2^n - 1)} x^n.
\]
By \eqref{eq:exp_gen_function_sol_1}, the r.h.s. of the above equation is equal to $F(x/2)$ and using \eqref{eq:expansion_of_Fr} we get 
\begin{equation}\label{eq:remainder}
R(x):=E(x)-\tilde{E}(x)=(e^x-1)\displaystyle \sum_{k=1}^\infty h \big( \frac12 (x/ {2^k}) \big)=:\sum_{k=1}^\infty R_k(x).
\end{equation}
Let us obtain a series expansion of $R(x)=\sum_{n=2}^\infty\frac{r_n}{n!}x^n$. 
Before going further, we first heuristically show the bound $r_n\le c(5/4)^n$. It can be shown that $h(x)$ is convex and increasing on the region $x\ge 0$ with $h(0)=0$. Therefore, for all $k\ge 2$ and $x>0$, $h(x/2^k)\le (1/2^{k-2})h(x/2^2)$. Using this together with the inequality $h(x)\le ce^{x}$, we get (constants may change with no harm)
$
|E(x)-\tilde{E}(x)|\le  c e^{\frac54x},
$
and this proves the bound. For a rigorous proof, however, we need some more efforts. Let's come back to the expansion of $R(x)$ in \eqref{eq:remainder}. 
By using the formula \eqref{eq:function_tilde_E_n}, we get
\begin{eqnarray*}
R_k(x)&=&(e^x-1)h \big( \frac12 (x/{2^k}) \big)\\
&=&\Big(\left(e^{x/2^k}\right)^{2^k}-1\Big)h \big( \frac12 (x/{2^k}) \big)\\
&=&\Big(1+\sum_{l=1}^{2^k-1}e^{lx/2^k}\Big)\left(e^{x/2^k}-1\right)h \big( \frac12 (x/{2^k}) \big)\\
&=&\Big(1+\sum_{l=1}^{2^k-1}e^{lx/2^k}\Big)
\frac13\sum_{s=2}^\infty\frac1{s!}\left(\frac32\right)^s(x/2^k)^s.
\end{eqnarray*}
Expanding the exponential function in the first term, we get
\begin{eqnarray*}
R_k(x)&=&\frac13\sum_{n=2}^\infty\frac1{n!}\left(\frac32\right)^n
\left(\frac1{2^k}\right)^nx^n
+\frac13\sum_{l=1}^{2^k-1}\sum_{m=0}^\infty\frac1{m!}\left(\frac{l}{2^k}\right)^mx^m\sum_{s=2}^\infty\frac1{s!}\left(\frac32\frac1{2^k}
\right)^sx^s\\
&=&\frac13\sum_{n=2}^\infty\frac1{n!}\left(\frac32\right)^n
\left(\frac1{2^k}\right)^nx^n
+\frac13\sum_{l=1}^{2^k-1}\sum_{n=2}^\infty \frac1{n!}\sum_{s=2}^n {n\choose s}\left(\frac32\right)^sl^{n-s}\left(\frac1{2^k}\right)^nx^n
\end{eqnarray*}
Therefore, we have
\begin{equation}\label{eq:remainder_coefficients}
r_n=\frac13\sum_{k=1}^\infty\left(\left(\frac32\right)^n
\left(\frac1{2^k}\right)^n+\sum_{l=1}^{2^k-1}\sum_{s=2}^n {n\choose s}\left(\frac32\right)^sl^{n-s}\left(\frac1{2^k}\right)^n\right).
\end{equation}
We exchange the order in the second summation: 
\[
\sum_{k=1}^\infty\sum_{l=1}^{2^k-1}\cdots=\sum_{m=0}^\infty
\sum_{l=2^m}^{2^{m+1}-1}\sum_{k=m+1}^\infty\cdots. 
\]
Then, summing over $k$ we get
\begin{equation}\label{eq:remainder_coefficients1}
r_n=\frac13\frac1{2^n-1}\left(\left(\frac32\right)^n+\sum_{s=2}^n{n\choose s}\left(\frac32\right)^s\sum_{m=0}^\infty
\sum_{l=2^m}^{2^{m+1}-1} l^{n-s}2^{-mn} \right).
\end{equation}
In the second term, we change the oder of summation: since $2^m\le l<2^{m+1}$, we have $\log_2l-1<m\le \log_2l$, or 
$m=[ \log_2l]$, where $[ a]$ is the integer part of $a$, i.e.,
\[
\sum_{m=0}^\infty
\sum_{l=2^m}^{2^{m+1}-1}\cdots=\sum_{l=1}^\infty\sum_{m=[ \log_2l]}^{[ \log_2l]}\cdots.
\]
Putting $[ \log_2l]=\log_2l-\delta(l)$ with $0\le \delta(l)<1$, we have $2^{-[ \log_2l]n}=l^{-n}2^{\delta(l)n}$.
Therefore, we have 
\begin{equation}\label{eq:remainder_coefficients2}
r_n=\frac13\frac1{2^n-1}\left(\left(\frac32\right)^n+\sum_{s=2}^n
{n\choose s}\left(\frac32\right)^s\sum_{l=1}^\infty l^{-s}2^{\delta(l)n}\right).
\end{equation}
This is the formula of the remainder in the statement of the theorem. We promptly see that $r_n$ is finite. 
In order to get the bound of $r_n$, let $N$ be a fixed large number which will be determined later. Let us consider the summation over $m$ in \eqref{eq:remainder_coefficients1}: $\sum_{m=0}^\infty a_m$, where
\[
a_m:=\sum_{l=2^m}^{2^{m+1}-1} l^{n-s}2^{-mn}.
\]
We see that 
\begin{eqnarray*}
a_{m+1}&=&\sum_{l=2^{m+1}}^{2^{m+2}-1} l^{n-s}2^{-(m+1)n}\\
&=&\sum_{l=2^m}^{2^{m+1}-1}\Big((2l)^{n-s}2^{-mn}2^{-n}+(2(l+1/2)^{n-s}2^{-mn}2^{-n}\Big)\\
&\le&\left(1+(1+\frac1{2^{m+1}})^{n-s}\right)2^{-s}a_m\\
&\le&\left(1+\frac1{2^{m+1}}\right)^{n}2^{-(s-1)}a_m\le e^{\frac{n}{2^{m+1}}}2^{-(s-1)}a_m.
\end{eqnarray*}
Therefore, for any $k\ge 1$, 
\[
a_{N+k}\le \left(\prod_{u=1}^ke^{\frac{n}{2^{N+u}}}\right)2^{-k(s-1)}a_N\le \left(e^{1/2^N}\right)^n2^{-k(s-1)}a_N\le (1+\e)^n2^{-k(s-1)}a_N,
\]
where we have taken $N$ large enough so that $e^{1/2^N}\le 1+\e$. Thus the second term inside the bracket in \eqref{eq:remainder_coefficients1} can be bounded by 
\begin{equation}\label{eq:r_n_upper_bound}
\sum_{s=2}^n{n\choose s}\left(\frac32\right)^s\left(1+\sum_{m=1}^{N}
\sum_{l=2^m}^{2^{m+1}-1} l^{n-s}2^{-mn}+(1+\e)^n \sum_{l=2^N}^{2^{N+1}-1} l^{n-s}2^{-Nn}\right).
\end{equation}
Now for any $1\le m\le N$, 
\begin{eqnarray}\label{eq:bound_each_term}
&&\sum_{s=2}^n{n\choose s}\left(\frac32\right)^s
\sum_{l=2^m}^{2^{m+1}-1} l^{n-s}2^{-mn}\nonumber\\
&\le &\sum_{s=2}^n{n\choose s}\left(\frac32\right)^s2^m2^{(m+1)(n-s)}2^{-mn}\nonumber\\
&\le&2^m\sum_{s=2}^n{n\choose s}\left(\frac3{2^{m+1}}\right)^s2^{n-s}\nonumber\\
&\le& 2^m\left(2+\frac3{2^{m+1}}\right)^n\le 2^N(11/4)^n=\mathrm{o}(3^n).
\end{eqnarray}
By \eqref{eq:remainder_coefficients1}, \eqref{eq:r_n_upper_bound}, and \eqref{eq:bound_each_term}, we have $r_n=\mathrm{o}(3/2)^n$. The proof is completed. 
\EndProof
 
\vskip 1 true cm
\noindent {\bf Acknowledgements}\\
We are grateful to Professor Tomoyuki Shirai for helping us with an improvement in asymptotics. 
The research by H.J.Yoo was supported by Basic Science Research Program through the National
Research Foundation of Korea (NRF) funded by the Ministry of Education (NRF-2016R1D1A1B03936006).

\end{document}